\documentclass{kluwer}

\begin{document}

\begin{article}

\begin{opening}

\title{An algorithm to obtain global solutions of the \\double
confluent Heun equation}

\author{J. \surname{Abad}}
\author{F. J. \surname{G\'omez}}
\author{J. \surname{Sesma}}
\institute{Departamento de F\'{\i}sica Te\'orica, Universidad de Zaragoza, Spain}


\dedication{In memoriam Prof. Luigi Gatteschi}

\runningtitle{Double confluent Heun equation}
\runningauthor{Abad e.a.}


\begin{abstract}
A procedure is proposed to construct solutions of the double confluent
Heun equation with a determinate behaviour at the singular points.
The connection factors are expressed as quotients of Wronskians
of the involved solutions. Asymptotic expansions are used in the
computation of those Wronskians. The feasibility of the method is shown
in an example, namely, the Schr\"odinger equation with a quasi-exactly-solvable
potential
\end{abstract}


\abbreviations{\abbrev{FDE}{Fuchsian differential equations};
  \abbrev{DCHE}{double confluent Heun equation}}


\classification{Mathematics Subject Classifications (2000)} {34B30,
33E20, 34M35}

\end{opening}

\section{Introduction}
Fuchsian differential equations (FDE), i. e., homogeneous
differential equations with single-valued analytic coefficient
functions, play a relevant role in many areas of Physics, Engineering
and Statistics. A large class of the special functions appearing
in these areas are but solutions
of second order FDE. The most studied of these equations is the Riemann
hypergeometric equation, that, in its general form, presents three
regular singular points. According with the location of its singularities,
the equation is known with different names (generalized Legendre equation,
Gauss equation, ...). Confluence of two or all three of its singularities
transforms the hypergeometric equation in others, with particular names
(Kummer, Bessel, Weber, ...), encountered in the solution of a lot of problems.

The next, in order of complexity, FDE is the Heun equation, that presents four
regular singularities. Different kinds of confluence of these singularities
produce the confluent, double confluent, biconfluent and triconfluent Heun equations.
The book edited by  \citeauthor{ronv} \shortcite{ronv} contains a thorough revision of the
results and open problems concerning those equations. More recently, \citeauthor{slav}
\shortcite{slav} have considered again the class of Heun equations, in their study of
Special Functions as solutions of differential equations with singular points. Confluence
reduces the number
of singularities, but increases the rank of some ones. Trivial changes of variable
allow to locate the singularities at chosen points. We report in Table \ref{singpts}
the type of singularities and their positions when those equations are written in
their usual conventional form .

\begin{table}
\caption[]{Singularities of the different types of Heun equations}
\label{singpts}
\begin{tabular}{lcccc}
\hline
equation & rank 0 & rank 1 & rank 2 & rank 3 \\
\hline
Heun & 0, 1, $a$ (arbitrary $\neq 0, 1, \infty$), $\infty$ &  &   &   \\
confluent H. &   0, 1 & $\infty$ &  &   \\
double confluent H. &      & 0, $\infty$ &   &   \\
biconfluent H. &   0 &   & $\infty$  &   \\
triconfluent H. &     &   &   & $\infty$ \\
\hline
\end{tabular}
\end{table}

In almost all nontrivial applications, in which special functions arise, the problem at hand
reduces to a boundary value one on an interval of the real axis whose interior points
are ordinary ones, whereas one or both of its ends
may be singular points of the differential equation. Usually, all but one of the parameters
appearing in the coefficient functions of the equation are fixed. Then, the boundary
conditions restrict (quantize) the possible
values of the free parameter. In the neighbourhood  of  each end, a pair of independent
solutions with a well known behaviour can be chosen. Such pair can serve as a basis to write,
by linear combination, any other solution of the differential equation, in particular, that
satisfying the boundary condition at the corresponding end. In order to examine the behaviour
of this solution at the other end, one has to be able to write each one of the two solutions
forming the basis at one end as linear combination of the basic solutions at the other end. In other
words, one has to solve the so called ``connection problem" for the ends. Articles by \citeauthor{kohn}
\shortcite{kohn}, \citeauthor{naun} \shortcite{naun}, \citeauthor{schs} \shortcite{schs}, and
\citeauthor{sch0} (\citeyear{sch0}, \citeyear{sch4}), among others,
deal with this problem.

In a recent paper \cite{gome} we have suggested a different procedure to solve the connection
problem for the cases of one of the ends being an ordinary or a regular singular point, and the other
an irregular one. Examples of such cases are found in the triconfluent and biconfluent Heun equations,
respectively. Our method is mostly inspired by the work of Naundorf, but presents some considerable
advantages, from the numerical point of view. Our use of the Wronskians of the basic solutions to
calculate the connection factors allows us to avoid having to sum slowly convergent double series like
those encountered in the Naundorf's procedure. The purpose of this paper is to extend our method
to the solution of the connection problem for two irregular singular points, the DCHE being a
typical example of this case.

We start by presenting, in the next Section, the DCHE and its basic solutions in the vicinity of
its singular points. Our procedure to solve the connection problem and the corresponding algorithms
are reported in Sections 3 and 4, respectively. An example, in Section 5, illustrates the method. 
Section 6 explains how to test the results of the only nontrivial algorithm. Finally, Section 7 
contains a comment about the generality of the procedure.

\section{The double confluent Heun Equation}
A complete discussion of the different forms adopted by the DCHE, under appropriate
transformations, has been done by \citeauthor{schw} \shortcite{schw}. We concentrate on what
they call {\em normal form} (Eq. 1.1.10 of their article)
\begin{equation}
D^2\,y + B(z)\,y = 0, \quad D=z\frac{d}{dz}, \quad B(z)=\sum_{p=-2}^2B_p\,z^p,  \label{uno}
\end{equation}
that, with the change of dependent variable
\begin{equation}
w(z)=z^{1/2}\,y(z)  \label{dos}
\end{equation}
and the notation
\begin{equation}
A_0=B_0+1/4,\qquad A_p=B_p, \quad p\neq 0, \label{tres}
\end{equation}
can be written in the form
\begin{equation}
z^2\,\frac{d^2w}{dz^2} + \sum_{p=-2}^2A_p\,z^p\, w=0, \label{cuatro}
\end{equation}
free of first order derivatives and, therefore, preferable for our discussion
because of the fact that the Wronskians of their solutions are independent of $z$.
In what follows, we will assume that we are dealing with a
{\em non-degenerate case} of the DCHE, that is, that
\[
A_2\,A_{-2} \neq 0.
\]
The origin and the infinity are the only singularities, both of rank 1, of the
differential equation.
We are interested in considering three pairs of independent solutions, namely:
\begin{itemize}
\item
Two {\em multiplicative} solutions \cite{arsc}, $w_1$ and $w_2$, that, except for particular sets of
values of the parameters $A_p$ in Eq. (\ref{cuatro}), have the form
\begin{equation}
w_j=z^{\nu_j}\sum_{n=-\infty}^{\infty} c_{n,j}\,z^n, \quad \mbox{being} \;
 \sum_{n=-\infty}^{\infty} |c_{n,j}|^2<\infty, \quad j=1,2. \label{cinco}
 \end{equation}
 The indices $\nu_j$ are not completely defined. They admit addition of any integer (with an adequate relabeling
 of the coefficients). To avoid ambiguities, we assume that
 \[
 |\Re\, \nu_j|\leq 1/2.
 \]
\item
 Two formal solutions, $w_3$ and $w_4$, that have the nature of asymptotic expansions for
 $z\to\infty$,
 \begin{equation}
w_k(z)\equiv\exp\left(\alpha_{k}\, z\right) z^{\mu_k}\,\sum_{m=0}^\infty
a_{m,k}\,z^{-m}, \quad a_{0,k}\neq 0, \quad k=3,4. \label{seis}
\end{equation}
It is usual to say that these two expansions are associated to each other.
\item
Two formal solutions, $w_5$ and $w_6$, asymptotic expansions for $z\to 0$, of the form
 \begin{equation}
w_l(z)\equiv\exp\left(\beta_{l}\, z^{-1}\right) z^{\rho_l}\,\sum_{m=0}^\infty
b_{m,l}\,z^{m}, \quad b_{0,l}\neq 0, \quad l=5,6. \label{siete}
\end{equation}
Also these expansions are associated.
\end{itemize}
The determinaton of the index $\nu$ and the coefficients $c_n$ of the
multiplicative solutions is rather laborious. By substitution of (\ref{cinco}) in (\ref{cuatro})
one obtains the infinite set of homogeneous equations for the coefficients
\begin{equation}
(n\! +\! \nu)(n\! -\! 1\! +\! \nu)\,c_{n} + \sum_{p=-2}^2A_p\,c_{n-p}=0,\quad n= \ldots, -1, 0, 1, \ldots,
 \label{ocho}
\end{equation}
that can be interpreted as a nonlinear eigenvalue problem, where the eigenvalue $\nu$ must
be such that
\begin{equation}
\sum_{n=-\infty}^{\infty} |c_{n}|^2<\infty. \label{nueve}
\end{equation}
In  Section 4 we recall the Newton iterative method to solve that problem. In general, two
indices, $\nu_1$ and $\nu_2$, and two corresponding sets of coefficients, $\{c_{n,1}\}$ and
$\{c_{n,2}\}$ are obtained, but for certain sets of values of the parameters $A_p$ only
one multiplicative solution appears. Any other independent solution must include powers of the variable
multiplied by its logarithm. Such logarithmic solutions cannot correspond, usually, to the
practical system that one tries to describe and are, therefore, to be discarded. We will assume,
from now on, that the parameters $A_p$ are such that the DCHE admits two independent multiplicative
solutions. In what concerns the formal solutions, the
exponents $\alpha$, $\mu$, $\beta$, and $\rho$ and the coefficients $a_m$ and $b_m$
of the expansions (\ref{seis}) and (\ref{siete}) must be such that these expansions
satisfy the differential equation. One obtains in this way
\begin{eqnarray}
\alpha_{k} & = & \sqrt{-A_{2}},\qquad \;\mu_k = -\,A_1/2\alpha_k,  \label{diez} \\
\beta_{l}  & = & \sqrt{-A_{-2}}, \qquad  \rho_l =  1 + A_{-1}/2\beta_l,  \label{duno}
\end{eqnarray}
for the exponents, and the recurrence relations
\begin{eqnarray}
\hspace{-.5cm}2\alpha_k\,m\,a_{m,k} & = & (m\! -\! \mu_k)(m\! -\! 1\! -\! \mu_k)\,a_{m-1,k}+\sum_{p=0}^2A_{-p}\,a_{m-1-p,k},
\label{ddos}  \\
2\beta_l\,m\,b_{m,l} & = & (m\! -\! 1\! +\! \rho_l)(m\! -\! 2\! +\! \rho_l)\,b_{m-1,l}+\sum_{p=0}^2A_{p}\,b_{m-1-p,l},
\label{dtres}
\end{eqnarray}
for the coefficients. The two independent solutions labeled by $k=3,4$ and those by
$l=5,6$ correspond to the two opposite values of the square root function in
(\ref{diez}) and (\ref{duno}).

\section{The connection factors}
Any solution $w$ of the DCHE can be written as a linear combination of the two multiplicative solutions,
\begin{equation}
w=\zeta_1\,w_1+\zeta_2\,w_2.  \label{dcuatro}
\end{equation}
Its behavior in the neighbourhood of the singular points can be immediately written if, besides the coefficients
$\zeta_1$ and $\zeta_2$, one knows the behaviour of the multiplicative solutions, that is, if one knows the
{\em connection factors} $T$ of their asymptotic expansions,
\begin{eqnarray}
w_j & \sim & T_{j,3}\,w_3+T_{j,4}\,w_4,\quad \mbox{for} \quad z\to\infty, \quad j=1, 2, \label{dcinco} \\
w_j & \sim & T_{j,5}\,w_5+T_{j,6}\,w_6,\quad \mbox{for} \quad z\to 0, \quad j=1, 2. \label{dseis}
\end{eqnarray}
These connection factors are, obviously, numerical constants, but their values
depend on the sector of the complex plane where $z$ lies. This fact, known as ``Stokes phenomenon"
\cite{ding}, introduces a slight complication in the procedure. As it is well known,
the connection factor multiplying each one of the asymptotic expansions in the right hand sides
of (\ref{dcinco}) and (\ref{dseis}) takes different values in the sectors of the complex $z$-plane
separated by a Stokes ray of the associated expansion. On the ray, the value of the connection factor
is the average of those two different ones. In the DCHE, both singular points have rank 1 and,
consequently, for each one of the expansions (\ref{seis}) and (\ref{siete}),
the sectors delimited by two contiguous Stokes rays of the same expansion have amplitude $2\pi$.
In fact, the Stokes rays for those expansions are given by
\begin{eqnarray}
\arg z & = & -\,\arg \alpha_k\pm 2n\pi,\quad n=0, 1, 2, \ldots,
\quad \mbox{for}\; w_k, \: k=3, 4, \label{dsiete}  \\
\arg z & = & \arg \beta_l\pm 2n\pi,\quad n=0, 1, 2, \ldots,
\quad \mbox{for}\; w_l, \; l=5, 6 \label{docho}
\end{eqnarray}
To be specific, let us assume that the labels $\{3, 4\}$ and $\{5, 6\}$ are
assigned in such a manner that
\begin{eqnarray}
\mbox{either}\quad -\pi\leq\arg \alpha_3 & < & -\pi/2 \quad  \mbox{or} \quad  \pi/2\leq\arg \alpha_3 < \pi,   \label{dnueve}  \\
-\pi/2\leq\arg \alpha_4 & < & \pi/2,  \label{veinte}  \\
\mbox{either}\quad -\pi<\arg \beta_5 & \leq & -\pi/2 \quad  \mbox{or} \quad \pi/2<\arg \beta_5 \leq \pi,  \label{vuno} \\
-\pi/2<\arg \beta_6 & \leq & \pi/2.  \label{vdos}
\end{eqnarray}
Then in the principal Riemann sheet, $-\pi<\arg z\leq\pi$, $T_{j,3}$ and  $T_{j,4}$
change their values as $z$ crosses respectively the rays $\arg z =-\, \arg \alpha_4$ and
$\arg z = -\,\arg \alpha_3$, whereas  $T_{j,5}$ and $T_{j,6}$ change at $\arg z = \arg \beta_6$ and
$\arg z = \arg \beta_5$, respectively.

Now, we are going to present our procedure to calculate the connection factors.
Let $\{r,s\}$ denote either the pair  $\{3,4\}$  or  $\{5,6\}$.
From Eqs. (\ref{dcinco}) or (\ref{dseis}), it is immediate to obtain
\begin{equation}
T_{j,r} = \frac{\mathcal{W}[w_j,w_s]}{\mathcal{W}[w_r,w_s]}, \qquad
T_{j,s} = \frac{\mathcal{W}[w_j,w_r]}{\mathcal{W}[w_s,w_r]},  \label{vtres}
\end{equation}
where $\mathcal{W}[f,g]$ represents the Wronskian of the functions $f$ and $g$,
\begin{equation}
\mathcal{W}[f,g](z) \equiv f(z)\,g^{\prime}(z) - f^{\prime}(z)\,g(z).   \label{vcuatro}
\end{equation}
There is no difficulty in using the asymptotic expansions (\ref{seis}) and (\ref{siete})
to calculate the denominators in (\ref{vtres}). In fact,
\begin{eqnarray}
\mathcal{W}[w_3,w_4] & = & -\, \mathcal{W}[w_4,w_3] = 2\,\alpha_4\,a_{0,3}\,a_{0,4},  \label{vcinco} \\
\mathcal{W}[w_5,w_6] & = & -\, \mathcal{W}[w_6,w_5] = 2\,\beta_5\,b_{0,5}\,b_{0,6}.  \label{vseis}
\end{eqnarray}
The numerators, instead, require a more sophisticated method. Direct evaluation would give
asymptotic expansions, containing both positive and negative powers of the variable, which
result to be inadequate to obtain precise values of those numerators. But we can benefit
from the fact that, whenever one considers not merely a ray, but a sector of the complex plane,
an asymptotic power series defines unambiguously an analytic function. Our idea is to write,
for each one of the connection factors $T_{j,t}$ ($t=r, s$), two
asymptotic expansions of the same function, one of the expansions containing,
as a common multiplicative constant, the Wronskian in the numerator of the expression of
$T_{j,t}$. Comparison of coefficients of equal powers of the variable in the two expansions
gives immediately that numerator. With that purpose, we introduce some auxiliary functions
and asymptotic expansions, namely
\begin{eqnarray}
\hspace{-.5cm} u_j(z) & \equiv & \mbox{e}^{-\alpha_tz/2}\,w_j,\quad u_t(z) \equiv \mbox{e}^{-\alpha_tz/2}\,w_t,
\quad j=1, 2, \quad t=3, 4, \label{vsiete}  \\
v_j(z) & \equiv & \mbox{e}^{-\beta_t/2z}\,w_j,\quad v_t(z) \equiv \mbox{e}^{-\beta_t/2z}\,w_t,
\quad j=1, 2, \quad t=5, 6, \label{vocho}
\end{eqnarray}
whose Wronskians obey the relations
\begin{eqnarray}
\mathcal{W}[u_j,u_t] & = & \mbox{e}^{-\alpha_tz}\, \mathcal{W}[w_j,w_t],\quad j=1, 2, \quad t=3, 4,   \label{vnueve} \\
\mathcal{W}[v_j,v_t] & = & \mbox{e}^{-\beta_tz^{-1}}\, \mathcal{W}[w_j,w_t], \quad j=1, 2, \quad t=5, 6.  \label{treinta}
\end{eqnarray}
Direct calculation of the left hand sides of these equations gives the  expansions
\begin{eqnarray}
\mathcal{W}[u_j,u_t] & = & \sum_{n=-\infty}^{\infty}\gamma_{n,j,t}\,z^{n+\nu_j+\mu_t},\quad j=1, 2, \quad t=3, 4,  \label{tuno} \\
\mathcal{W}[v_j,v_t] & = & \sum_{n=-\infty}^{\infty}\eta_{n,j,t}\,z^{n+\nu_j+\rho_t},\quad j=1, 2, \quad t=5, 6,  \label{tdos}
\end{eqnarray}
with coefficients
\begin{eqnarray}
\gamma_{n,j,t} & = &\!\!\sum_{m=0}^\infty a_{m,t}[\alpha_t\,c_{n+m,j}
-(n\! +\! 2m\! +\! 1\! +\! \nu_j\! -\! \mu_t)c_{n+m+1,j}],  \label{ttres} \\
\eta_{n,j,t} &  = & \!\!\sum_{m=0}^\infty b_{m,t}[-\beta_t\,c_{n-m+2,j}-
(n\! -\! 2m\! +\! 1\! +\! \nu_j\! -\! \rho_t)c_{n-m+1,j}],  \label{tcuatro}
\end{eqnarray}
Now, in order to compare coefficients of equal powers in the two sides, we need to write expansions
of the right hand sides  of (\ref{vnueve}) and (\ref{treinta}) in powers of $z$ with the same exponents
as in the left hand sides. This can be done by using the Heaviside's exponential series \cite{hard}
\begin{equation}
\mbox{e}^\xi \sim \sum_{n=-\infty}^\infty \frac{\xi^{n+\delta}}{\Gamma (n+1+\delta)}, \qquad |\arg \xi |<\pi,
\quad \delta \; \;\mbox{arbitrary},   \label{tcinco}
\end{equation}
with $\xi$ replaced respectively by $-\alpha_tz$ and $-\beta_tz^{-1}$, and taking $\delta$ equal to
$\nu_j+\mu_k$ and $-(\nu_j+\rho_t)$. This kind of representation of the exponential function
has already been used by Naundorf in his
solution of the connection problem. Notice that the restriction $|\arg \xi |<\pi$ prevents the use of
such expansions if $\arg z=-\arg \alpha_t$ or $\arg z=\arg \beta_t$. But these values of $\arg z$
correspond precisely to the rays at which the connection factor to be calculated changes.
As we have already said, the value assigned to the connection factor in that situation is the
average of its values in the sectors separated by that ray, where the restriction $|\arg \xi |<\pi$
is fulfilled. Following that procedure, one obtains, for $z$ out of the Stokes ray of $w_t$
in the principal Riemann sheet,
\begin{eqnarray}
\mathcal{W}[w_j,w_t] & = & \frac{\Gamma (n\! +\! 1+\! \nu_j\! +\! \mu_t)}{(-\alpha_t)^{n+\nu_j+\mu_t}}\,\gamma_{n,j,t},
\quad j=1, 2, \quad t=3, 4,  \label{tseis} \\
\mathcal{W}[w_j,w_t] & = & \frac{\Gamma (n\! +\! 1\! -\! \nu_j\! -\! \rho_t)}{(-\beta_t)^{n-\nu_j-\rho_t}}\,\eta_{-n,j,t},
\quad j=1, 2, \quad t=5, 6,  \label{tsiete}
\end{eqnarray}
where the minus sign in front of $\alpha_t$ and $\beta_t$ is to be interpreted as $\mbox{e}^{i\pi}$ or
$\mbox{e}^{-i\pi}$ so as to have $|\arg (-\alpha_tz)|<\pi$ and $|\arg (-\beta_tz^{-1})|<\pi$. If $z$ lies
on the Stokes ray of $w_t$, one has, instead of (\ref{tseis}) and (\ref{tsiete}),
\begin{eqnarray}
\mathcal{W}[w_j,w_t] & = & (-1)^n\,\cos[\pi(\nu_j+\mu_t)]\,
\frac{\Gamma (n\! +\! 1+\! \nu_j\! +\! \mu_t)}{\alpha_t^{n+\nu_j+\mu_t}}\,\gamma_{n,j,t},
 \label{tocho} \\
\mathcal{W}[w_j,w_t] & = & (-1)^n\,\cos[\pi(\nu_j+\rho_t)]\,
\frac{\Gamma (n\! +\! 1\! -\! \nu_j\! -\! \rho_t)}{\beta_t^{n-\nu_j-\rho_t}}\,\eta_{-n,j,t}.
  \label{tnueve}
\end{eqnarray}

\section{The algorithms}
We present in this Section algorithms that facilitate implementation of the procedure
sketched above. Its successive steps  are considered separately in different subsections

\subsection{Multiplicative solutions}
As it has been already said, we face in Eqs. (\ref{ocho}) and (\ref{nueve}) a nonlinear eigenvalue problem.
Algorithms to solve finite order problems of this kind have been discussed by \citeauthor{ruhe} \shortcite{ruhe}.
The condition (\ref{nueve}) implies
\begin{equation}
\lim_{n\to\pm\infty}\,|c_n|= 0,  \label{cuarenta}
\end{equation}
that makes possible to reduce, by truncation, our problem (\ref{ocho}) to one with
$n$ going from $-M$ to $N$, both $M$ and $N$ being positive and sufficiently large to ensure
that the solution of the truncated problem approximates that of the original one.
Each step of the Newton iteration method consists in moving from an approximate solution,
$\{\nu^{(i)}, c_n^{(i)}\}$, to another one,  $\{\nu^{(i+1)}, c_n^{(i+1)}\}$, by solving the
system of equations
\begin{eqnarray}
\hspace{-.5cm}\Big(2n\! -\! 1\! +\! 2\nu^{(i)}\Big)c_n^{(i)}\Big(\nu^{(i+1)} &\! -&\!  \nu^{(i)}\Big) +
\Big(n\! +\! \nu^{(i)}\Big)\Big(n\! -\! 1\! +\! \nu^{(i)}\Big)\,c_n^{(i+1)}
\nonumber \\ +\,\sum_{p=-2}^2A_p\,c_{n-p}^{(i+1)} & = & 0, \qquad n=-M, \ldots, -1, 0, 1, \ldots, N,
\label{cuno} \\
\sum_{n=-M}^N{c_n^{(i)}}^*c_n^{(i+1)} & = & 1,  \label{cdos}
\end{eqnarray}
that results, by linearization, from (\ref{ocho}) and from the truncated normalization condition
\[
\sum_{n=-M}^N\left| c_{n}\right|^2 = 1.
\]
Needless to say,
the values of $c_m^{(i)}$ with $m<-M$ or $m>N$ entering in some of the equations (\ref{cuno}) should be taken
equal to zero, in accordance with the truncation done. As usual, the iteration process is stopped when the
difference between consecutive solutions is satisfactory. Then the process is repeated with
larger and larger values of $M$ and $N$, to obtain a stable solution within the required precision.

The iteration process just described needs initial values $\{\nu^{(0)}, c_n^{(0)}\}$ not far from the true solution.
The two different values of $\nu$ can be obtained from the two eigenvalues
\begin{equation}
\lambda_j = \exp (2i\pi\nu_j)  \label{ctres}
\end{equation}
of the circuit matrix $\mathbb{C}$ \cite{waso} for the singular point at $z=0$.
The entries of that matrix can be computed by numerically integrating the
DCHE on the unit circle, from $z=\exp (0)$ to $z=\exp (2i\pi)$, for two independent sets of
initial values. If we consider two solutions, $w_a(z)$ and $w_b(z)$, obeying, for instance,
the conditions
\begin{eqnarray}
w_a(\mbox{e}^0)&=&1,\qquad w_a^\prime(\mbox{e}^0)=0, \nonumber \\
w_b(\mbox{e}^0)&=&0,\qquad w_b^\prime(\mbox{e}^0)=1, \nonumber
\end{eqnarray}
then
\begin{eqnarray}
C_{11}&=&w_a(\mbox{e}^{2i\pi}), \qquad C_{12}=w_b(\mbox{e}^{2i\pi}), \nonumber \\
C_{21}&=&w_a^\prime(\mbox{e}^{2i\pi}), \qquad C_{22}=w_b^\prime(\mbox{e}^{2i\pi}), \nonumber
\end{eqnarray}
and
\begin{equation}
\nu=\frac{1}{2i\pi}\,\ln \left[\frac{1}{2}\left( C_{11}+C_{22}\pm\sqrt{\left(C_{11}\! -\! C_{22}\right)^2+4C_{12}C_{21}}\right)\right].
\label{ccuatro}
\end{equation}
The two signs in front of the square root produce two different values for $\nu$, unless the parameters
$A_p$ in the DCHE be such that $\left(C_{11}-C_{22}\right)^2+4C_{12}C_{21}=0$, in which case only one
multiplicative solution appears, any other independent solution containing logarithmic terms.
The ambiguity in the real part of $\nu$ due to the multivaluedness of the logarithm in the right
hand side of (\ref{ccuatro}) is eliminated by the restriction $|\Re\, \nu|\leq 1/2$ assumed above.
Notice that
\[
\lambda_1\,\lambda_2=\det \mathbb{C}=\mathcal{W}[w_a,w_b]=1
\]
and, therefore,
\[
\nu_1+\nu_2=0\quad (\mbox{mod} \;1).
\]
This may serve as a test for the integration of the DCHE on the unit circle.

Although Eq. (\ref{ccuatro}) is exact, the $C_{mn}$ are obtained numerically and the resulting values of
$\nu$ may only be considered as starting values, $\nu_j^{(0)}$ ($j=1, 2$), for the Newton iteration process.
As starting coefficients $c_n^{(0)}$ one may use the solutions of the homogeneous system
\begin{eqnarray}
\lefteqn{\hspace{-1cm}(n\! +\! \nu_j^{(0)})(n\! -\! 1\! +\! \nu_j^{(0)})\,c_{n,j}^{(0)} +
\sum_{p=-2}^2A_p\,c_{n-p,j}^{(0)}=0}& & \nonumber \\
& & \hspace{2cm}n= -M, \ldots, -1, 0, 1, \ldots, N, \qquad j=1, 2,  \label{ccinco}
\end{eqnarray}
with the already mentioned truncated normalization condition
\begin{equation}
\sum_{n=-M}^N |c_{n,j}^{(0)}|^2 = 1.  \label{cseis}
\end{equation}

\subsection{Formal solutions}
The exponents of the formal expansions (\ref{seis}) and (\ref{siete}) are given in (\ref{diez})
and (\ref{duno}). The coefficients must obey the recurrence relations (\ref{ddos}) and
(\ref{dtres}) that are but third order difference equations. The Perron-Kreuser theorem
\cite{perr} predicts for each one of them a unique (save for multiplication by a constant) dominant
solution that can be obtained starting, for instance, with
\begin{equation}
a_{0,k}=1, \quad b_{0,l}=1, \qquad k=3,4, \quad l=5,6,  \label{csiete}
\end{equation}
and using the recurrence relations as they appear in (\ref{ddos}) and (\ref{dtres}).
To avoid overflows, it may be convenient to deal with the quotients of successive coefficients, that satisfy
\begin{eqnarray}
2\alpha_k\,m\,\textstyle{\frac{a_{m,k}}{a_{m-1,k}}} & = & (m\! -\! \mu_k)(m\! -\! 1\! -\! \mu_k)+A_0  \nonumber \\
& & \hspace{2cm}+\,\frac{A_{-1}}{\frac{a_{m-1,k}}{a_{m-2,k}}}+\frac{A_{-2}}{\frac{a_{m-1,k}}{a_{m-2,k}}\,\frac{a_{m-2,k}}{a_{m-3,k}}},
\label{cocho}  \\
2\beta_l\,m\,\textstyle{\frac{b_{m,l}}{b_{m-1,l}}} & = & (m\!-\! 1\!  +\! \rho_l)(m\! -\! 2\! +\! \rho_l)+A_0 \nonumber \\
& & \hspace{2cm}+\,\frac{A_{1}}{\frac{b_{m-1,l}}{b_{m-2,l}}}+\frac{A_{2}}{\frac{b_{m-1,l}}{b_{m-2,l}}\,\frac{b_{m-2,l}}{b_{m-3,l}}},
\label{cnueve}
\end{eqnarray}

\subsection{Connection factors}
For the computation of the connection factors one should make use of Eqs. (\ref{vtres}), (\ref{vcinco}), (\ref{vseis}),
(\ref{tseis}) and (\ref{tsiete}), complemented by (\ref{ttres}) and  (\ref{tcuatro}) conveniently truncated.
The integer $n$ in Eqs. (\ref{ttres}),  (\ref{tcuatro}),  (\ref{tseis})  and
(\ref{tsiete}) may be chosen at will, with the only restriction of being positive and satisfying
\begin{equation}
(n+\nu)(n+\nu-1) > \sum_{p=-2}^2 |A_p|.   \label{cincuenta}
\end{equation}
Use of different values of $n$ may serve as a test of the procedure. The sums in (\ref{ttres}) and (\ref{tcuatro}), truncated and
written in terms of quotients of successive coefficients, $a_{m,k}/a_{m-1,k}$ and $b_{m,l}/b_{m-1,l}$, may be computed in nested form.

\section{An example}
A particular case of DCHE on the positive real semiaxis, $z\in [0,
+\infty )$, is the Schr\"odinger equation
\begin{equation}
-\frac{\hbar^2}{2m}\,\left(\frac{d^2
R(r)}{dr^2}-\frac{l(l+1)}{r^2}\,\,R(r)\right)+V(r)\,R(r)=E\,R(r),
\label{extrauno}
\end{equation}
for the reduced radial wave function $R(r)$ of
a particle of mass $m$, angular momentum $l\hbar$ and energy
$E=A_2\hbar^2/2mr_0^2$ in a spherically symmetric potential
\begin{equation}
V(r)=-\,\frac{\hbar^2}{2m}
\left(\frac{A_{-2}\,r_0^2}{r^4}+\frac{A_{-1}\,r_0}{r^3}+
\frac{A_{0}+l(l\! +\! 1)}{r^2}+\frac{A_{1}\,r_0^{-1}}{r}\right).
\label{quno}
\end{equation}
In fact, by using the variable $z$ and the wave function $w$ given
by
\begin{equation}
z=r/r_0  \qquad \mbox{and} \qquad w(z)=R(r), \label{extrados}
\end{equation}
the Schr\"odinger equation (\ref{extrauno}) adopts the form (\ref{cuatro}).
The potential (\ref{quno}), with some restrictions on the values of the parameters,
belongs to a class of quasi-exactly-solvable ones
\cite{turb}. In fact, for the particular set of parameters
\begin{equation}
A_{-2}= -\, 1, \quad A_{-1}= 4/5, \quad A_{0}=31/25, \quad
A_{1}=3/5, \label{qdos}
\end{equation}
it presents an $l=0$ bound state of energy $E=-\,(1/4)\,\hbar^2/2mr_0^2$ \cite{ozce}. In other words,
Eq. (\ref{cuatro}), with the values of the $A_p$ given by (\ref{qdos}) and by $A_2=-1/4$,
has a normalizable solution on the positive real semiaxis. To illustrate the procedure described
in the preceding Sections,
we have applied it to find global solutions, for $z\in [0,+\infty)$, of (\ref{cuatro}) with fixed
values (\ref{qdos}) of the
potential parameters and several values
\begin{equation}
A_2 = -\,1/10, \, -\,1/5, \,-\,1/4, \,-\,3/10, \,-\,2/5,
\label{qtres}
\end{equation}
of the energy parameter. The results are shown in Tables \ref{multi} to \ref{conne2}. We report the output of our
double precision FORTRAN codes, but, due to roundoff errors, we do not claim that all the digits reproduced are
correct. In fact, entries that in the tables appear as having modulus less than $10^{-12}$ should be exactly
equal to zero. Table \ref{multi} shows the index $\nu_1$ of the multiplicative solution $w_1$. Of course, the
index $\nu_2$ of $w_2$ is the opposite, $\nu_2=-\nu_1$. The connection factors of both multiplicative solutions
with the formal ones are listed in Table \ref{conne1}. These factors depend on the normalization adopted
for the different solutions. We have taken
\[
c_{0,j}=1, \qquad j=1, 2,
\]
for the multiplicative solutions, whereas
the normalization of the formal solutions is determined by
\[
a_{0,k}=1, \quad k=3, 4, \qquad b_{0,l}=1, \quad  l=5, 6.
\]
Table \ref{regu} gives the coefficients $\zeta_1$ and $\zeta_2$ of a linear combination of the multiplicative
solutions,
\[
w_{\mbox{reg}} = \zeta_1\,w_1 + \zeta_2\,w_2,
\]
well behaved near the origin on the positive real semiaxis, that is, such that
\[
w_{\mbox{reg}}(z) \sim w_5(z) \qquad \mbox{as}\; z\to 0^+.
\]
Finally, Table \ref{conne2} shows the connection factors
giving the behaviour of $w_{\mbox{reg}}$ at infinity on the positive real semiaxis,
\[
w_{\mbox{reg}}(z) \sim T_{\mbox{reg},3}\,w_3(z) + T_{\mbox{reg},4}\,w_4(z) \qquad \mbox{as}\; z\to +\infty.
\]
\begin{table}
\caption[]{Indices of the multiplicative solutions of the DCHE for the sets of parameters given in (\ref{qdos}) and
(\ref{qtres}). The two indices have opposite values; so the real and imaginary parts of only one of them are shown.}
\label{multi}
\begin{tabular}{lrr}
\hline
\multicolumn{1}{c}{$A_2$} & \multicolumn{1}{c}{\hspace{1cm}$\Re\, \nu_1$} & \multicolumn{1}{c}{\hspace{0.75cm}$\Im\, \nu_1$}  \\
\hline
$-$1/10 & \hspace{2cm}$-$.500000000000E+00 & \hspace{1.5cm}$-$.703150555392E+00  \\
$-$1/5 & $-$.500000000000E+00 & $-$.531738153810E+00\\
$-$1/4 & $-$.400000000000E+00 & .102867462041E$-$26 \\
$-$3/10 & $-$.262974969075E$-$15 & .509507933497E+00 \\
$-$2/5 & .181198729462E$-$16 & .688682990633E+00 \\
\hline
\end{tabular}
\end{table}
\begin{table}
\caption[]{Connection factors of the multiplicative solutions, $w_1$ and $w_2$, with the formal solutions
at infinity, $w_3$ and $w_4$, and at the origin, $w_5$ and $w_6$, for the DCHE (\ref{cuatro})
with parameters (\ref{qdos}) and (\ref{qtres}).}
\label{conne1}
\begin{tabular}{lrrrr}
\hline
\multicolumn{1}{c}{$A_2$} & \multicolumn{1}{c}{$\Re\, T_{1,3}$} & \multicolumn{1}{c}{$\Im\, T_{1,3}$}
& \multicolumn{1}{c}{$\Re\, T_{1,4}$} & \multicolumn{1}{c}{$\Im\, T_{1,4}$}  \\
\hline
{\tiny $-$0.1 }& {\tiny $-$.130166356702E+01} & {\tiny .218375080422E+01}
& {\tiny .537554017411E+00} & {\tiny .314875723054E+00}\\
{\tiny $-$0.2} & {\tiny $-$.204239771333E+01} & {\tiny .362437608911E+00}
& {\tiny .462340626214E$-$01} & {\tiny .190707221682E+00} \\
{\tiny $-$0.25} & {\tiny $-$.786334477859E+00} & {\tiny .770822392694E$-$13}
& {\tiny .227716282126E$-$13} & {\tiny .938023374790E$-$28} \\
{\tiny $-$0.3} & {\tiny .107546538975E+00} & {\tiny .137759967302E+01}
 & {\tiny $-$.106218738591E+00} & {\tiny .108713582942E$-$01} \\
{\tiny $-$0.4} & {\tiny .108556195947E+01} & {\tiny .199397670471E+01}
& {\tiny $-$.125046997516E+00} & {\tiny .673923636677E$-$01} \\
\hline
\multicolumn{1}{c}{$A_2$} & \multicolumn{1}{c}{$\Re\, T_{2,3}$} & \multicolumn{1}{c}{$\Im\, T_{2,3}$}
& \multicolumn{1}{c}{$\Re\, T_{2,4}$} & \multicolumn{1}{c}{$\Im\, T_{2,4}$}  \\
\hline
{\tiny $-$0.1 }& {\tiny $-$.179473952568E+00} & {\tiny .141771950373E+01}
& {\tiny $-$.347068293144E+00} & {\tiny $-$.466342783505E$-$01}\\
{\tiny $-$0.2} & {\tiny .104106884981E+01} & {\tiny .818532573406E+00}
& {\tiny $-$.834077925914E$-$01} & {\tiny .934810009512E$-$01} \\
{\tiny $-$0.25} & {\tiny .491548162537E+00} & {\tiny .488371898613E$-$13}
& {\tiny $-$.289740699463E$-$01} & {\tiny $-$.293358750439E$-$14} \\
{\tiny $-$0.3} & {\tiny .107546538975E+00} & {\tiny $-$.137759967302E+01}
 & {\tiny $-$.106218738591E+00} & {\tiny $-$.108713582942E$-$01} \\
{\tiny $-$0.4} & {\tiny .108556195947E+01} & {\tiny $-$.199397670471E+01}
& {\tiny $-$.125046997516E+00} & {\tiny $-$.673923636677E$-$01} \\
\hline
\multicolumn{1}{c}{$A_2$} & \multicolumn{1}{c}{$\Re\, T_{1,5}$} & \multicolumn{1}{c}{$\Im\, T_{1,5}$}
& \multicolumn{1}{c}{$\Re\, T_{1,6}$} & \multicolumn{1}{c}{$\Im\, T_{1,6}$}  \\
\hline
{\tiny $-$0.1 }& {\tiny $-$.350619618820E+01} & {\tiny $-$.780802661791E+00}
& {\tiny .283770014388E$-$01} & {\tiny $-$.136520514662E+00}\\
{\tiny $-$0.2} & {\tiny $-$.218278236288E+01} & {\tiny $-$.183368690581E+00}
& {\tiny .349930523633E$-$02} & {\tiny $-$.829409172034E$-$01} \\
{\tiny $-$0.25} & {\tiny $-$.786334477859E+00} & {\tiny $-$.770774921662E$-$13}
& {\tiny .581790855708E$-$15} & {\tiny .371252705024E$-$27} \\
{\tiny $-$0.3} & {\tiny $-$.185903592530E+00} & {\tiny .135138118813E+01}
 & {\tiny .590575874067E$-$01} & {\tiny .525872360505E$-$02} \\
{\tiny $-$0.4} & {\tiny $-$.379273027228E$-$01} & {\tiny .219921373757E+01}
& {\tiny .927424524732E$-$01} & {\tiny .159396234117E$-$03} \\
\hline
\multicolumn{1}{c}{$A_2$} & \multicolumn{1}{c}{$\Re\, T_{2,5}$} & \multicolumn{1}{c}{$\Im\, T_{2,5}$}
& \multicolumn{1}{c}{$\Re\, T_{2,6}$} & \multicolumn{1}{c}{$\Im\, T_{2,6}$}  \\
\hline
{\tiny $-$0.1 }& {\tiny .182308938715E+01} & {\tiny .867938208517E+00}
& {\tiny .346916923655E$-$01} & {\tiny $-$.702844924334E$-$01}\\
{\tiny $-$0.2} & {\tiny .128432288745E+01} & {\tiny .553414307656E+00}
& {\tiny .229815763899E$-$01} & {\tiny $-$.477580284950E$-$01} \\
{\tiny $-$0.25} & {\tiny .493413857062E+00} & {\tiny $-$.477076860664E$-$13}
& {\tiny .144870349731E$-$01} & {\tiny $-$.142012480382E$-$14} \\
{\tiny $-$0.3} & {\tiny $-$.185903592530E+00} & {\tiny $-$.135138118813E+01}
 & {\tiny .590575874067E$-$01} & {\tiny $-$.525872360505E$-$02} \\
{\tiny $-$0.4} & {\tiny $-$.379273027228E$-$01} & {\tiny $-$.219921373757E+01}
& {\tiny .927424524732E$-$01} & {\tiny $-$.159396234117E$-$03} \\
\hline
\end{tabular}
\end{table}
\begin{table}
\caption[]{Coefficients of the linear combination of multiplicative solutions resulting in a solution
 $w_{\mbox{reg}}=\zeta_1w_1+\zeta_2w_2$ normalizable, on the positive real semiaxis, at $z=0$.
 The parameters of the DCHE are the same as in Tables \ref{multi} and \ref{conne1}.}
\label{regu}
\begin{tabular}{lrrrr}
\hline
\multicolumn{1}{c}{$A_2$} & \multicolumn{1}{c}{$\Re\, \zeta_1$} & \multicolumn{1}{c}{$\Im\, \zeta_1$}
& \multicolumn{1}{c}{$\Re\, \zeta_2$} & \multicolumn{1}{c}{$\Im\, \zeta_2$}  \\
\hline
{\tiny $-$0.1 }& {\tiny $-$.136295750328E+00} & {\tiny .283302821758E$-$01}
& {\tiny .222074702056E+00} & {\tiny $-$.109613756594E+00}\\
{\tiny $-$0.2} & {\tiny $-$.228256434722E+00} & {\tiny .963021587152E$-$02}
& {\tiny .322449452131E+00} & {\tiny $-$.155165465358E+00} \\
{\tiny $-$0.25} & {\tiny $-$.127172345631E+01} & {\tiny .124655928872E$-$12}
& {\tiny .510716705831E$-$13} & {\tiny .325901912099E$-$25} \\
{\tiny $-$0.3} & {\tiny $-$.325468718708E$-$01} & {\tiny $-$.365514500226E+00}
 & {\tiny $-$.325468718708E$-$01} & {\tiny .365514500226E+00} \\
{\tiny $-$0.4} & {\tiny $-$.390741170002E$-$03} & {\tiny $-$.227347243108E+00}
& {\tiny $-$.390741170002E$-$03} & {\tiny .227347243108E+00} \\
\hline
\end{tabular}
\end{table}
\begin{table}
\caption[]{Connection factors of the ``regular" solution
at the origin with the formal solutions $w_3$ and $w_4$ for $z\to\infty$ along
the positive real semiaxis. The parameters in the DCHE are the same as in the preceding Tables.}
\label{conne2}
\begin{tabular}{lrrrr}
\hline
\multicolumn{1}{c}{$A_2$} & \multicolumn{1}{c}{$\Re\, T_{reg,3}$} & \multicolumn{1}{c}{$\Im\, T_{reg,3}$}
& \multicolumn{1}{c}{$\Re\, T_{reg,4}$} & \multicolumn{1}{c}{$\Im\, T_{reg,4}$}  \\
\hline
{\tiny $-$0.1 }& {\tiny .231089872113E+00} & {\tiny $-$.485000928307E$-$12}
& {\tiny $-$.164373692458E+00} & {\tiny .317038062470E$-$13}\\
{\tiny $-$0.2} & {\tiny .925400135830E+00} & {\tiny .898559004980E$-$12}
& {\tiny $-$.247795480195E$-$01} & {\tiny .560940183192E$-$13} \\
{\tiny $-$0.25} & {\tiny .100000000000E+01} & {\tiny $-$.196048546485E$-$12}
& {\tiny $-$.304389678921E$-$13} & {\tiny .162523413961E$-$26} \\
{\tiny $-$0.3} & {\tiny .100006470514E+01} & {\tiny .000000000000E+00}
 & {\tiny .148614535378E$-$01} & {\tiny .000000000000E+00} \\
{\tiny $-$0.4} & {\tiny .905801865774E+00} & {\tiny .000000000000E+00}
& {\tiny .307406581930E$-$01} & {\tiny .000000000000E+00} \\
\hline
\end{tabular}
\end{table}

We have already mentioned that the DCHE (\ref{cuatro}) with
parameters as given by (\ref{qdos}) and by $A_2=-1/4$ possesses a
solution,
\begin{equation}
w(z) = z^{3/5}\,\exp(-z^{-1}-z/2), \label{qcuatro}
\end{equation}
normalizable on the positive real semiaxis. This fact allows to know the exact values of
the connection factors $T_{1,3}$ and $T_{1,5}$ of the multiplicative solution
\begin{equation}
w_1(z) = \sum_{n=-\infty}^\infty c_{n,1}\,z^{n-2/5}.  \label{qcinco}
\end{equation}
Writing (\ref{qcuatro}) in the form
\begin{equation}
w(z) = \sum_{n=-\infty}^\infty \hat{c}_{n}\,z^{n+3/5}, \label{qseis}
\end{equation}
with
\[
\hat{c}_n = \left\{ \begin{array}{ll}(-1)^n\displaystyle{\sum_{m=0}^\infty\frac{2^{-m}}{m!(m-n)!}}
 \quad \mbox{for}\; n<0, \\
(-1)^n\displaystyle{\sum_{m=0}^\infty\frac{2^{-m-n}}{m!(m+n)!}}
\quad \mbox{for}\; n\geq 0, \end{array} \right.
\]
comparison with (\ref{qcinco}) gives immediately
\[
T_{1,3} = T_{1,5} = c_{0,1}/\hat{c}_{-1}
\]
and, having chosen $c_{0,1}=1$,
\[
T_{1,3} = T_{1,5} = -\,\left[\sum_{m=0}^\infty \frac{2^{-m}}{m!(m+1)!}\right]^{-1}.
\]
The values of $\Re\, T_{1,3}$ and $\Re\, T_{1,5}$ obtained with our procedure and reported in
Table \ref{conne1} coincide with the exact values up to all shown digits.

\section{A test of the multiplicative solutions}

We have described, in subsection 4.1, an algorithm to obtain the
multiplicative solutions, $w_1$ and $w_2$, mentioned in Eq.
(\ref{cinco}). A test of the correct implementation of the algorithm
could be to try to reproduce the particular solution of the DCHE
provided by the Maple 11 system. The statement
\[
\mbox {HeunD}(\alpha, \beta, \gamma, \delta, t)
\]
gives the value at $t$ ($|t|<1$) of a function $y(t)$ that obeys
the DCHE, written in the Jaff\'{e}-Lay form \cite{slav},
\begin{equation}
\frac{d^2y(t)}{dt^2}-\frac{\alpha+2t+\alpha t^2-2t^3}{(t^2-1)^2}\,\,
\frac{dy}{dt} + \frac{\delta+(2\alpha+\gamma)t+\beta
t^2}{(t^2-1)^3}\,\,y = 0,  \label{extratres}
\end{equation}
with singular points located at $t=-1$ and $t=1$, and satisfies the
boundary conditions
\begin{equation}
y(0)=1, \qquad y^\prime (0)=0.  \label{extracuatro}
\end{equation}
(We adopt the notation of the Maple manual. Needless to say, the parameters represented 
by the symbols $\alpha$ and $\beta$ along this Section are not related to the exponents 
$\alpha_k$, $k=3, 4$, and $\beta_l$, $l=5, 6$ of the formal solutions (\ref{seis}) and 
(\ref{siete}).) The changes of independent and dependent variables
\begin{equation}
t=\frac{z-1}{z+1} \quad \mbox{and} \quad
y(t)=z^{-1/2}\exp\left(\frac{\alpha}{8}\left(z-\frac{1}{z}\right)\right)\,w(z)
\label{extracinco}
\end{equation}
transform Eq. (\ref{extratres}) in Eq. (\ref{cuatro}) with
parameters
\[
A_{-2}=-\frac{\alpha^2}{64}, \quad
A_{-1}=\frac{\gamma-\beta-\delta}{16}, \quad
A_{0}=\frac{8-\alpha^2+4(\beta-\delta)}{32}, \nonumber
\]
\begin{equation}
A_{1} = \frac{-\gamma-\beta-\delta}{16}, \quad
A_{2}=-\frac{\alpha^2}{64}.  \label{extraseis}
\end{equation}
Obviously, $w(z)$ in (\ref{extracinco}) can be written as a linear
combination
\begin{equation}
w(z)=\xi_1\,w_1(z) + \xi_2\,w_2(z)   \label{extraextra}
\end{equation}
of the two multiplicative solutions, the coefficients $\xi_1$ and
$\xi_2$ depending on the parameters $\alpha$, $\beta$, $\gamma$ and
$\delta$. Then, the procedure to check the algorithm used to obtain
the multiplicative solutions would run along the following steps:
\begin{itemize}
\item
Choose a set of values for the parameters $\alpha$, $\beta$,
$\gamma$ and $\delta$.
\item
Compute the corresponding values of $A_p$, $p=-2,-1,\ldots,2$, by
means of the relations (\ref{extraseis}).
\item
Use the algorithm described in Subsection 4.1 to obtain the two
multiplicative solutions.
\item
Determine the coefficients $\xi_1$ and $\xi_2$ in (\ref{extraextra})
by requiring that the second of Eqs. (\ref{extracinco}) and that
resulting by deriving it with respect to $t$ be satisfied for $t=0$
and $z=1$. In other words, solve the system
\begin{eqnarray}
w_1(z=1)\,\xi_1 + w_2(z=1)\,\xi_2 & = & 1,  \label{extrasiete}   \\
w_1^\prime(z=1)\,\xi_1 + w_2^\prime(z=1)\,\xi_2 & = & -\, 1/2.
\label{extraocho}
\end{eqnarray}
\item
Choose an arbitrary value of $t$ inside the unit circle and compute
the corresponding value of $z$, according to the first of Eqs.
(\ref{extracinco}).
\item
Check the fulfilment of the second of Eqs. (\ref{extracinco}) for
those values of $t$ and $z$, the left hand side being computed by
the above mentioned Maple statement and the right hand side by using
the multiplicative solutions and the coefficients $\xi_1$ and
$\xi_2$.
\end{itemize}
Let us choose, for instance, the parameters
\[
\alpha=4, \quad \beta=-\,3, \quad \gamma=2, \quad \delta=-\,1,
\]
for which the corresponding parameters in Eq. (\ref{cuatro}) are
\[
A_{-2}=-1/4, \quad A_{-1}=3/8, \quad A_{0}=-1/2, \quad A_{1} = 1/8,
\quad A_{2}=-1/4,
\]
The multiplicative solutions, obtained with a double precision
FORTRAN code, are then those given in Table \ref{ultima}, and the
coefficients in (\ref{extraextra}) giving the $w(z)$ corresponding
to the Maple particular solution turn out to be
\[
\xi_1=-\,0.48013092979925, \qquad \xi_2=1.5087428376316.
\]
\begin{table}
\caption[]{Coefficients $c_{n,j}$ of the multiplicative solutions,
$w_j(z)=z^{\nu_j}\sum_{n=-\infty}^{\infty} c_{n,j}\,z^n$, $j=1,2$,
normalized in such a way that $c_{0,1}=c_{0,2}=1$. Only the most
relevant coefficients are shown. The indices are
$\nu_1=-\,\nu_2=0.346120772343$.} \label{ultima}
\begin{tabular}{lrrrr}
\hline \multicolumn{1}{c}{$m$} & \multicolumn{1}{c}{$c_{-m,1}$} &
\multicolumn{1}{c}{$c_{m,1}$} & \multicolumn{1}{c}{$c_{-m,2}$} &
\multicolumn{1}{c}{$c_{m,2}$} \\
\hline
1 & {\tiny .167014439267E+00} & {\tiny .190993398030E+01} & {\tiny $-$.107279310787E+00} & {\tiny .382083714867E+00} \\
2 & {\tiny .481375140054E$-$01} & {\tiny $-$.501405710995E$-$02} & {\tiny .395518603946E$-$01} & {\tiny .343615502071E+00} \\
3 & {\tiny .256845652795E$-$02} & {\tiny    .651050357942E$-$01} & {\tiny $-$.297089364185E$-$02} & {\tiny .126405802915E$-$01} \\
4 & {\tiny .670742341737E$-$03} & {\tiny $-$.688158578777E$-$03} & {\tiny    .484111513331E$-$03} & {\tiny .916616846207E$-$02} \\
5 & {\tiny .151109427662E$-$04} & {\tiny    .719872336102E$-$03} & {\tiny $-$.276617916112E$-$04} & {\tiny .120052304745E$-$03} \\
6 & {\tiny .436467584453E$-$05} & {\tiny $-$.788297814632E$-$05} & {\tiny    .284956022379E$-$05} & {\tiny .881920361011E$-$04} \\
7 & {\tiny .424159499646E$-$07} & {\tiny    .392398781019E$-$05} & {\tiny $-$.131313020620E$-$06} & {\tiny .507159044452E$-$06} \\
8 & {\tiny .163574045519E$-$07} & {\tiny $-$.405522033871E$-$07} & {\tiny    .982771078139E$-$08} & {\tiny .435958940659E$-$06} \\
9 & {\tiny .537668556535E$-$10} & {\tiny    .127234781325E$-$07} & {\tiny $-$.379605967044E$-$09} & {\tiny .109230687414E$-$08} \\
10& {\tiny .397572715764E$-$10} & {\tiny $-$.122029907678E$-$09} & {\tiny    .222381579709E$-$10} & {\tiny .131080043658E$-$08} \\
11& {\tiny $-$.119341203300E$-$13} & {\tiny    .273442588328E$-$10} & {\tiny $-$.739683873557E$-$12} & {\tiny .105746798402E$-$11} \\
12& {\tiny    .676604265157E$-$13} & {\tiny $-$.243166143611E$-$12} & {\tiny    .355327272482E$-$13} & {\tiny .264897017572E$-$11} \\
13& {\tiny $-$.164661027834E$-$15} & {\tiny    .417999894189E$-$13} & {\tiny $-$.103815479335E$-$14} & {\tiny $-$.464027978928E$-$15} \\
14& {\tiny    .850622285001E$-$16} & {\tiny $-$.345793838238E$-$15} & {\tiny    .422141078690E$-$16} & {\tiny    .384446086066E$-$14} \\
15& {\tiny $-$.319253539778E$-$18} & {\tiny    .477715538676E$-$16} & {\tiny $-$.109995830505E$-$17} & {\tiny $-$.299696460688E$-$17} \\
16& {\tiny    .821887767466E$-$19} & {\tiny $-$.369227215875E$-$18} & {\tiny    .387438635716E$-$19} & {\tiny    .420067427815E$-$17} \\
17& {\tiny $-$.376967703210E$-$21} & {\tiny    .423582849020E$-$19} & {\tiny $-$.911209342015E$-$21} & {\tiny $-$.490269012114E$-$20} \\
18& {\tiny    .629190033709E$-$22} & {\tiny $-$.307215713853E$-$21} & {\tiny    .282929069113E$-$22} & {\tiny    .358011576237E$-$20} \\
19& {\tiny $-$.321864802770E$-$24} & {\tiny    .299866441211E$-$22} & {\tiny $-$.606469439157E$-$24} & {\tiny $-$.509135489559E$-$23} \\
20& {\tiny    .390955943808E$-$25} & {\tiny $-$.204922402639E$-$24} & {\tiny    .168292193160E$-$25} & {\tiny    .244636867112E$-$23} \\
\hline
\end{tabular}
\end{table}

\section{Final comment}

The procedure presented in the preceding Sections referred to global
solutions in the interval between two adjacent irregular singular
points, of rank 1, of a second order differential equation like the
DCHE. It can be trivially extended, however, to the case of two
irregular singular points of arbitrary rank. The only changes due to
the larger rank of the singularities are in the evaluation of the
Wronskians of the pairs formed by each one of the multiplicative
solutions and each one of the formal solutions. The steps are analogous 
to those detailed in a paper \cite{gome} dealing with
the connection problem for an ordinary or regular singular point and
an irregular one of arbitrary rank. Obviously, the
series solutions around the ordinary or regular singular point 
in that reference should be replaced by the 
multiplicative solutions of the problem at hand. Then, having computed the connection
factors of each multiplicative solution with the formal solutions at
each singular point, the procedure is the same as described above.

\section*{Acknowledgments}

Our thanks are due to an anonymous referee of a previous version of
this paper for valuable comments contributing to improve
considerably the presentation of our algorithm and for suggesting us
to test our procedure in the form described in Section 6. We
acknowledge the financial support of  Comisi\'{o}n
Inter\-mi\-nis\-te\-rial de Cien\-cia y Tec\-no\-lo\-g\'{\i}a and of
Di\-pu\-ta\-ci\'on Ge\-ne\-ral de Ara\-g\'on

\end{article}

\end{document}